\documentclass[11pt]{amsart}
\usepackage{amsmath,amsfonts,latexsym,amssymb,amscd, graphicx,pictexwd}

\newcommand{\concat}{\ensuremath{+\!\!\!\!+\,}}
\renewcommand{\P}{{\mathbb P}}
\newcommand{\E}{{\mathbb E}}
\newcommand{\R}{{\mathbb R}}

\newcommand{\ZZ}{{\mathbb Z}}

\newcommand{\F}{{\mathbb F}}
\newcommand{\G}{{\mathbb G}}

\newcommand{\HH}{{\mathbf H}}
\newcommand{\0}{{\mathbf 0}}

\newcommand{\cc}{{\mathbf c}}
\newcommand{\zz}{{\mathbf z}}

\newcommand{\yy}{{\mathbf y}}
\newcommand{\mm}{{\mathbf m}}

\newcommand{\xx}{{\mathbf x}}
\newcommand{\nn}{{\mathbf n}}
\newcommand{\ee}{{\mathbf e}}
\newcommand{\dd}{{\mathbf d}}

\newcommand{\cal}{\mathcal}

\newcommand{\calR}{{\mathcal R}}

\newcommand{\calZ}{{\mathcal Z}}

\newcommand{\Exp}{{\rm Exp}}

\newcommand{\argmax}{\mathop{\rm arg\,max}}

\newtheorem{thm}{Theorem}
\newtheorem{lem}{Lemma}

\numberwithin{equation}{section}

\begin{document}
\title[Duality between coalescence times and exit points]{Duality between coalescence times and exit points in last-passage percolation models}

\author{Leandro P. R. Pimentel}
\address{Universidade Federal do Rio de Janeiro\\
Caixa Postal 68530, 21941-909 Rio de Janeiro, RJ, Brasil}
\email{leandro@im.ufrj.br}

\thanks{The author was partially supported by the CNPQ grant 474233/2012-0.}

\date{\today}
%\AMSsubject{}
%\keywords{}

\begin{abstract}
 In this article we prove a duality relation between coalescence times and exit points in last-passage percolation models with exponential weights. As a consequence, we get lower bounds for coalescence times, with scaling  exponent $3/2$, and we relate its distribution with variational problems involving the Brownian motion process and the Airy$_2$ process. The proof relies on the relation between Busemann functions and the Burke property for stationary versions of the last-passage percolation model with boundary.        
\end{abstract}

\maketitle

\section{introduction and main results}

\subsection{Introduction}
This article develops and studies the scaling behavior of coalescence times of directional geodesics in the two-dimensional last-passage (site) percolation model with exponential passage times. In this context, geodesics are paths which locally maximizes the passage time, and a directional geodesic is a semi-infinite (oriented) geodesic that has an asymptotic direction. Uniqueness and coalescence of directional geodesics initially appeared in the scene of (first-passage) percolation in the work of Newmann and coauthors \cite{Ne}. These semi-infinite paths are the building blocks of Busemann functions and have become a key notion in the study of the geometry of percolation and equilibrium measures of related particle systems \cite{BaCaSe,CaPi,Co,FeMaPi}.  

A fundamental question concerns the exact scaling behavior of coalescence times. This random variable is conjectured to have scaling exponent $3/2$ and, under this scaling, it is also expected to converge to an universal limiting distribution. This conjecture is motivated by universality aspects of last-passage percolation models \cite{KZ}, which are known to lie in the Kardar-Parisi-Zhang universality class \cite{Jo1}. Although this question is of great interest, in the last two decades, the only available result in the literature was due to W\"utrich \cite{Wu}, where it was proved that the scaling exponent is greater than $3/2-\epsilon$ for all small $\epsilon>0$.     

The main contribution of this article resides in to bring new rigorous results on the subject, and also to shed new lights on the scaling scenario of coalescence times. We extend W\"utrich's result up to the $m^{3/2}$ scale, where $m$ denotes the distance between the starting points of the geodesics, by proving that the probability for coalescence after time $r m^{3/2}$ is of order $1-cr^2$, for small enough $r>0$ (Theorem \ref{thm:coaltail}). And we also prove that, for large enough $r>0$, this tail probability scales as fast as $r^{-2/3}$, which shows that the limiting coalescence time has heavy tail behavior (Theorem \ref{thm:lowtail}). In particular, one gets a non-integrable random variable.

The method of proof is based on a duality formula that relates the distribution of the coalescence time, with the probability that the ``exit-point counting measure'' of $[-m,m]$ is zero (Theorem \ref{thm:dualfor}). This counting measure is composed by exit points of a slightly different last-passage percolation model with stationary boundary conditions. This formula allows us to prove results for coalescence times by studying the respective problem in the dual context, where a considerable part of the fluctuation theory for exit-points is well developed \cite{BaCaSe,FlQuaRe,Jo}. Another interesting aspect of  duality is that it provides a description \eqref{conj1} of the conjectured scaling scene for coalescence times in terms of variational problems involving the Brownian motion and Airy$_2$ processes, which have received considerable attention recently \cite{QuaRe} from a more analytical point of view.

On our way to duality we derive some other results. For example, we prove \eqref{thm:Level}, which states that the level set of a Busemann function is distributed as a Palm version of the stationary totally asymmetric simple exclusion  process (TASEP) conditioned to have a jump at time zero. We also show self-duality of the directional geodesic tree (Lemma \ref{thm:selfdual}), which is a central property to derive the duality formula.     

We expect that duality will lead us to sharp upper bounds for coalescence times as well. However, there are still some technical obstacles for proving such bounds, and we will leave it for a future work. It is also natural to expect that the directional geodesics tree has a continuum scaling limit. This limit will not be the usual Brownian web, since they live in distinct universality classes. From physical motivation, and  connections with polymer models, we call this conjectured limit the (zero temperature) polymer web \eqref{polyweb}. Duality should also hold for the limiting polymer web.

\subsection*{Overview} In sections \ref{sec:dual} and \ref{sec:lower} we will formally introduce the LPP models and state the main results. In section \ref{sec:proofs} we will prove them. In section \ref{sec:final} we will make some more comments on the conjectured scaling scenario.

\subsection*{Acknowledgments}
The author would like to say thanks to Ivan Corwin and Timo Sepp\"al\"ainen for illuminating discussions and helpful comments on a previous version of this manuscript. 

\subsection{Duality between coalescence times and exit points}\label{sec:dual} 
Consider a collection of i.i.d. random variables $\omega=\{W_\xx\,:\,\xx\in\ZZ^2\}$, distributed according to an exponential distribution function of parameter one. In last-passage site percolation (LPP) models, each number $W_\xx$ represents the passage (or percolation) time through  vertex $\xx=(x(1),x(2))$. For $\ZZ^2$ lattice vertices $\xx\leq \yy$ (i.e. $x(i)\leq y(i)\,,i=1,2$), denote $\Gamma(\xx,\yy)$ the set of all up-right oriented paths $\gamma=(\xx_0,\xx_1\dots,\xx_k)$ from $\xx$ to $\yy$, i.e. $\xx_0=\xx$, $\xx_k=\yy$ and $\xx_{j+1}-\xx_j\in\{\ee_1,\ee_2\}$, for $j=0,\dots,k-1$, where $\ee_1:=(1,0)$ and $\ee_2=(0,1)$. The \emph{passage time} of $\gamma=(\xx_0,\xx_1\dots,\xx_k)\in\Gamma(\xx,\yy)$ is defined as 
$$W(\gamma):=\sum_{j=1}^k W_{\xx_i}\,.$$ 
The \emph{last-passage time} between $\xx$ and $\yy$ is defined as 
\begin{equation*}%\label{last}
L(\xx,\yy)=L_\omega(\xx,\yy):=\max_{\gamma\in\Gamma(\xx,\yy)}W(\gamma)\,.
\end{equation*}
(The lower index indicates the dependence with the environment $\omega$.) The \emph{geodesic} from $\xx$ to $\yy$ is the a.s. unique maximizing path $\gamma(\xx,\yy)=\gamma_\omega(\xx,\yy)\in\Gamma(\xx,\yy)$ such that   
$$L(\xx,\yy)=W(\gamma(\xx,\yy))\,.$$

The study of semi-infinite geodesics in last-passage percolation models with exponential weights was done in \cite{Co,FePi}. We summarize below the properties that we will use in this paper. The \emph{semi-infinite geodesic} starting at $\xx$ and along direction $\dd:=(1,1)$ is the almost surely  unique up-right oriented path $\gamma(\xx)=(\xx_n)_{n\geq 0}$ which satisfies:
\begin{itemize}
\item[(i)] $\xx_0=\xx$ and, for any $m<n$,
$$\gamma(\xx_m,\xx_n)=(\xx_{m},\cdots,\xx_n)\,;$$ 
\item[(ii)] For any sequence of lattice points $(\yy_n)_{n\geq 1}$ such that 
$$\yy_n=(y_n(1),y_n(2))\geq  \0=(0,0)\,\,\mbox{ and }\,\,\lim_{n\to\infty}\frac{y_n(2)}{y_n(1)}=1\,,$$ 
we have
$$\lim_{n\to\infty}\gamma(\xx,\xx+\yy_n)=\gamma(\xx)\,.$$
\end{itemize}
Convergence of a sequence of finite paths $(\gamma_n)_{n\geq 0}$ to a semi-infinite path $\gamma$ means that, for every finite set $K\subseteq\ZZ^2$, $\gamma_n$ and $\gamma$ will  coincide inside $K$, eventually. Another important property of semi-infinite geodesics with the same direction is coalescence. The symbol $\concat$ below stands for the concatenation of two paths.
\begin{itemize} 
\item[(iii)] For any $\xx,\yy\in\ZZ^2$ there exists $\cc\in\ZZ^2$ such that
$$\gamma(\xx)=\gamma(\xx,\cc)\concat\gamma(\cc)\,\,\mbox{ and }\,\, \gamma(\yy)=\gamma(\yy,\cc)\concat\gamma(\cc)\,.$$
\end{itemize}
We could have consider directional geodesics in a arbitrary fixed direction $\dd_a=(1,a)$, for $a\in(0,\infty)$ but, for the sake of simplicity, we will restrict our attention to $a=1$.

We note that if $\cc$ satisfies (iii), and $\cc'\in\gamma(\cc)$, then $\cc'$ also satisfies (iii). From now on we  denote $\cc(\xx,\yy)$ the first (in the up-right orientation) \emph{coalescence point}, in the sense that $\cc'\geq \cc(\xx,\yy)$ for every other geodesic point $\cc'$ that satisfies (iii). For $m\geq 1$, denote $\mm^h:=(m,0)$ and $\mm^v:=(0,m)$ and  let
$$T_m:=\mbox{the second coordinate of $\cc(\mm^h,\mm^v)$}\,.$$
By symmetry, it is clear that the first coordinate of $\cc_m$ has the same distribution as $T_m$. We call $T_m$ the \emph{coalescence time}. 
   
Now consider a slightly different LPP model, where we introduce boundary conditions as follows. Denote $\Exp(\rho)$ an exponential random variable with parameter $\rho$. Take an environment $\bar\omega=\{\bar W_{\zz}\,:\,\zz\geq 0\}$ mutually independent with the following distribution: 
$$\bar W_\zz\stackrel{dist.}{:=}\left\{\begin{array}{ll} 0\,, & \mbox{ if }\, \zz =\0\,;\\
\Exp(1)\,, & \mbox{ if }\, \zz > \0\,;\\
\Exp(1/2)\,,& \mbox{ otherwise }\,.\end{array}\right.$$
In other words, we put i.i.d. exponentials random variables of parameter $1/2$ along the horizontal and vertical axes of the first quadrant, and leave its interior with the same distribution as before. We denote 
$$\bar L(\xx):=L_{\bar\omega}(\0,\xx)\,$$
the last-passage time from $\0$ to $\xx$, with respect to the $\bar\omega$ environment. This LPP model can be seen as a stationary version of the classical LPP model previously introduced. As an effect of the boundary condition, we have that
$$\bar L(y,n)-\bar L(x,n)\stackrel{dist.}{=}\bar L(y,0)-\bar L(x,0)\stackrel{dist.}{=}\sum_{z=x+1}^y\bar W_{(z,0)}\,,$$
for all $n\geq 0$ and $x<y$.

We call the \emph{exit-point} of the geodesic $\gamma_{\bar\omega}(\0,\xx)$ (with respect to the environment $\bar\omega$) the last boundary point of the path (following the up-right orientation). To distinguish between exit via the horizontal or the vertical axis, we introduce a non-zero integer-valued random variable $Z(\xx)=Z_{\bar\omega}(\xx)$ such that if $Z(\xx)>0$ then the exit-point is $(Z(\xx),0)$, while if $Z(\xx)<0$ then the exit-point is $(0,-Z(\xx))$. The \emph{exit-point counting measure process} is defined as 
\begin{equation*}%\label{defpoint}
\cal Z_n:=(\zeta_n(z)\,,z\in\ZZ)\in\{0,1\}^{\ZZ}\,,
\end{equation*}  
where, for fixed $n\geq 1$,
$$\zeta_n(z):=\left\{\begin{array}{ll}1 & \mbox{ if $z=Z(x,n)$ for some $x\in [1,\infty)$}\,,\\ 
0 & \mbox{otherwise}\,.\end{array}\right.$$ 
The associated counting measure is defined as
$$\cal Z_n(A):=\sum_{z\in A} \zeta_n(z)\,,\,\mbox{ for }A\subseteq\ZZ \,.$$ 
 
The key result of this article is the following duality formula:
\begin{thm}\label{thm:dualfor} For $n>m>0$ we have that
\begin{equation}\label{dualfor}
\P\left(T_m< n\right)=\P\left(\cal Z_n\left([-m,m]\right)=0\right)\,.
\end{equation}
\end{thm}

\subsection*{A few words about the proof.} The coalescence property (iii) of semi-infinite geodesics allows us to introduce \emph{Busemann functions} in the LPP model, which are defined as 
\begin{equation*}%\label{Busemann}  
 B(\xx,\yy):=L(\yy,\cc)-L(\xx,\cc)\,,\mbox{ for }\xx,\yy\in\ZZ^2\,,
\end{equation*}
where $\cc=\cc(\xx,\yy)$. Busemann functions provide an alternative construction of stationary LPP models with boundary \cite{CaPi}. These models enjoys a very special property, named, the \emph{Burke property}. This property, formulated in terms of Busemann functions, will lead us to self-duality of the directional geodesic tree, composed by semi-infinite coalescent directional geodesics, and finally to the duality formula \eqref{dualfor}. 

\subsection{Lower bounds for the tail distribution}\label{sec:lower}
The tail distribution of the coalescence time is defined as\footnote{The reason to put the additional scaling factor $2^{-5/2}$ will became clear in the sequel, and it is related to universality of the expected limiting distribution.} 
$$\G(r):=\liminf_{m\to\infty}\P\left(\frac{T_m}{2^{-5/2}m^{3/2}}> r\right)\,.$$
We combine \eqref{dualfor} with scaling of exit points \cite{BaCaSe} to study the behavior of $\G$ close to $0$.
\begin{thm}\label{thm:coaltail}
There exist constants $c_0,r_0>0$ such that for all $r\in[0,r_0]$ we have 
\begin{equation}\label{coaltail}
\G(r)\geq 1-c_0 r^2\,.
\end{equation}
In particular,
\begin{equation*}
\lim_{r\to 0}\G(r)=1\,.
\end{equation*}
\end{thm}

The fluctuations of last-passage times are related to variational problems involving the Brownian motion and the Airy$_2$ process \cite{QuaRe}. The Airy$_{2}$ process is a one-dimensional stationary process with continuous paths, whose finite dimensional distributions are describe by Fredholm determinants. Duality allows us to link coalescence times with these processes as well. Let 
$$ U:=\argmax_{u\in\R}\left\{\sqrt{2} \cal B(u)+ \cal A(u)-u^2 \right\}\,,$$
where $(\cal B(u)\,,\,u\in\R)$ is a standard two-sided Brownian motion, and $(\cal A(u)\,,\,u\in\R)$ is an independent Airy$_2$ process, and denote
$$\F(s):=\P\left(U\leq s\right)\,.$$
The random variable $U$ is well defined (the location is a.s. unique \cite{Pi}), and it describes the limit in distribution of the rescaled exit point (Lemma \ref{thm:scalargmax}): 
$$\lim_{n\to\infty}\frac{Z(n,n)}{2^{5/3}n^{2/3}}\stackrel{dist.}{=}U\,.$$ 
Together with duality, this yields to: 
\begin{thm}\label{thm:lowtail} 
For $r>0$, we have that
\begin{equation}\label{lowtail}
\G(r)\geq \F(r^{-2/3})-\F(-r^{-2/3})\,.
\end{equation}
\end{thm}

A straightforward consequence of Theorem \ref{thm:lowtail} is that
\begin{equation*}%\label{eq:lowertail}
\liminf_{r\to\infty}r^{2/3}\G(r)\geq 2 f(0)\,,
\end{equation*}
where $f$ is the density of $\F$. Although, as far as the author knows, there is no analytical  description of $f$. We do expect that $f$ is bell shaped around $0$, as in the case of a Brownian motion minus a parabola \cite{Gro}, as well as in the case of an Airy$_2$ process minus a parabola \cite{FlQuaRe}. In particular, we also expect that  $f(0)>0$, which would imply non-integrability of $\G$. 

We conjecture that  \footnote{See Section \ref{sec:final} for further discussions on the conjectured picture for the scaling limit of coalescence times.}  
\begin{equation}\label{conj1}
\G(r)= \P\left(\cal U(-r^{-2/3},r^{-2/3}] \geq 1\right)\,,
\end{equation}
where $\cal U$ is a couting process composed by Dirac measures located at $U(v)$, for $v\in\R$, which is defined as
$$U(v):=\sup\argmax_{u\in\R}\{\sqrt{2} \cal B(u)+ \cal A(u,v)-(u-v)^2\}\,.$$
The process $(\cal A(u,v)\,,\,u,v\in\R$ is the so called Airy$_2$ sheet \cite{CoQua}. If this conjecture is true, then
\begin{equation}\label{conj2}
\lim_{r\to\infty}r^{2/3}\G(r)= \lim_{\delta\to 0}\delta^{-1}\P\left(\cal U(-\delta,\delta] \geq 1\right)\,.
\end{equation}
\newline

\section{Burke's property, level sets, self-duality and scaling of coalescence times}\label{sec:proofs}

\subsection{The last-passage percolation model and the exclusion process}
The LPP model can be seen as a function of the motion of particles in the one-dimensional \emph{totally asymmetric simple exclusion process} (TASEP). This process is a Markov process $(\eta_t\,,\,t\geq 0)$ in the state space $\{0,1\}^\ZZ$ whose elements are particle configurations: $\eta_t(x)=1$ indicates a particle at site $x$ at time $t$; otherwise $\eta_t(x)=0$ (a hole is at site $j$ at time $t$). With rate $1$, if there is a particle at site $x$, it attempts to jump to site $x+1$; if there is a hole at $x+1$ the jump occurs, otherwise nothing happens. The generator of the process is given by
$$\mathcal{G}f(\eta)=\sum_{x\in\ZZ}\eta(x)(1-\eta(x + 1))\left [ f(\eta^{x,x+1}) - f(\eta)\right]\,,$$
where $\eta^{x,y}(x)=\eta(z)$ $\forall z\not\in\{x,y\}$, $\eta^{x,y}(x)=\eta(y)$ and $\eta^{x,y}(y)=\eta(x)$.
For $p\in(0,1)$, let $\nu_p$ denote the product measure on $\ZZ$ with density $p$. Then $\nu_p$ is invariant for $\mathcal{G}$. The reverse process with respect to $\nu_p$ has generator $\mathcal{G}^*$ which is also a TASEP with reversed jumps:  
$$\mathcal{G}^*f(\eta)=\sum_{x\in\ZZ}\eta(x)(1-\eta(x-1))\left [ f(\eta^{x,x-1}) - f(\eta)\right]\,.$$
This property is called reversibility (or Burke's property) of the TASEP.

A construction of the (time) stationary process $\underline\eta=(\eta_t)_{t\in\R}$ with the marginal distribution $\nu_p$ can be done by choosing a configuration $\eta$ according to $\nu_p$ and then running the process with generator $L$ forward in time and the process with generator $L^*$ backward in time. The reversed process $\underline\eta^*$ is given by $\eta^*_t=\eta_{-t^-}$. The particle jumps of $\underline\eta$ induce a stationary point process $S$ in $\ZZ\times \R$. Let $S_x\subseteq\R$ be the (discrete and random) set of times for which a particle of $\underline\eta$ jumps from $x$ to $x+1$, and $S=(S_x\,,\,x\in\ZZ)$. The map $\underline\eta\mapsto S$ associates alternate point processes to each trajectory. The law of the process S is space and time translation invariant. Let $S^0$ be the Palm version of $S$, that is, the process with the law of $S$ conditioned to have a point at $(x, t)=(0, 0)$. In the corresponding process $\underline\eta^0=(\eta^0_t)_{t\in\R}$ there is a particle jumping from $0$ to $1$ at time zero. In the reverse process $\underline\eta^{*0}$ there is a particle jumping from $1$ to $0$ at time zero. 
 
We now construct a random function $\underline G=G(\underline\eta^0)$ as follows \cite{Ro}: first label the particles of $\eta^0_0$ in decreasing order, giving label $0$ to the particle at site $1$. We note that, for $\underline\eta^0$, at time zero the particle already has jumped from site $0$ to $1$. Call $P_j(0)$ the position of the $j$th particle at time zero; we have $P_0(0)=1$ and $P_{j+1}(0)<P_j (0)$ for all $j\in\ZZ$. Label the holes of $\eta^0_0$ in increasing order, giving the label $0$ to the hole at site $0$: $H_0(0) = 0$ and $H_{i+1}(0) > H_i(0)$ for all $i\in\ZZ$. The position of the $j$th particle and the $i$th hole at time t are denoted, respectively, $P_j(t)$ and $H_i(t)$. The order is preserved at later and earlier times: $P_j(t) > P_{j+1}(t)$ and $H_i(t) < H_{i+1}(t)$, for all $t\in\R$, $i,j\in\ZZ$. Let $G(i,j)$ denote the time the $i$th hole and the $j$th particle of $\underline\eta^0$ interchange positions; in particular $G(0, 0) = 0$. Let 
$$\underline G=G(\underline\eta^0):=\{G(\zz)\,:\, \zz\in\ZZ^2\}\,.$$ 
The LPP model with boundary condition and the TASEP are related by (we take $p=1/2$)
$$\{\bar L(\zz)\,:\, \zz\in\ZZ_+^2\}\stackrel{dist.}{=}\{G(\zz)\,:\, \zz\in\ZZ_+^2\}\,.$$  
For a proof of this distributional equality we refer to (4.21), (4.22) and Lemma 4.2 in Section 4.2 of \cite{FeMaPi}.

To construct the analog object for the reversed process, we set $\hat\eta^{*0}_t(j):=\eta^*_t(-j)$. By reversibility, $\hat{\underline\eta}^{*0}$ is also a stationary TASEP, but now with jumps in the same orientation as before. For this process, there is a particle jumping from $-1$ to $0$ at time zero. At this time, we give label $0$ to the particle at $0$ and label $0$ to hole at $-1$, and construct the interchanging times $G^*(i,j)$ as before, so that 
$$\underline G^{*}:=\{G^*(\zz)\,:\, \zz\in\ZZ^2\}=\{-G(-\zz)\,:\,\zz\in\ZZ^2\}\,.$$
As a consequence of reversibility $\underline G^{*}\stackrel{dist.}{=}\underline G$, and therefore
\begin{equation}\label{reverse}  
 \{-G(-\zz)\,:\,\zz\in\ZZ^2\}\stackrel{dist.}{=}\{G(\zz)\,:\, \zz\in\ZZ^2\}\,
\end{equation}

\subsection{Burke's property for Busemann functions}
In Cator and Pimentel \cite{CaPi}, it was developed a connection between the LPP model with boundary and Busemann functions. Almost sure existence and coalescence of semi-infinite geodesics along the negative diagonal direction are also true. Let $\gamma^{\downarrow}(\xx)=(\xx_n)_{n\geq 0}$ denote the down-left oriented semi-infinite geodesic starting at $\xx$ and  along the negative diagonal direction. Thus, $\gamma^{\downarrow}(\xx)$ satisfies (i), (ii) and (iii), but now in the down-left orientation. For $\xx,\yy\in\ZZ^2$, let $\cc^{\downarrow}(\xx,\yy)$ denote the coalescence point between $\gamma^{\downarrow}(\xx)$ and $\gamma^{\downarrow}(\yy)$, and set 
$$ B^{\downarrow}(\xx,\yy):=L(\cc^{\downarrow},\yy)-L(\cc^{\downarrow},\xx)\,.$$  
The main result in \cite{CaPi} states that
\begin{equation*}%\label{eqBu}
\{\bar L(\zz)\,:\,\zz\in\ZZ_+^2\}\stackrel{dist.}{=}\{B^{\downarrow}\left(\zz\right)\,:\,\zz\in \ZZ_+^2\}\,,
\end{equation*}
where $ B^{\downarrow}(\zz):=B^{\downarrow}\left(\0,\zz\right)$. It also follows from the results in \cite{CaPi} that
\begin{equation}\label{eqBu}
\{ G(\zz)\,:\,\zz\in\ZZ^2\}\stackrel{dist.}{=}\{ B^{\downarrow}\left(\zz\right)\,:\,\zz\in \ZZ^2\}\,.
\end{equation}
We call \eqref{BusBur} below the Burke property of Busemann functions. 

\begin{lem}\label{thm:BuBu} For the Busemman functions, we have that
\begin{equation}\label{BusBur}
\{B^{\downarrow}(\zz)\,:\,\zz\in\ZZ^2\}\stackrel{dist.}{=}\{- B^{\downarrow}(-\zz)\,:\,\zz\in\ZZ^2\}\,.
\end{equation}
\end{lem}

\begin{proof} It follows directly from \eqref{reverse} and \eqref{eqBu}.

\end{proof} 
  
\subsection{A remark on level sets of Busemann functions}\label{sec:Level}
Equation \eqref{eqBu}, which relates Busemann functions with the stationary TASEP conditioned to have a particle jumping from $0$ to $1$ at time zero, can be used to study level sets of Busemann functions (also called horospheres). Indeed, it is expected that the lattice boundary of the region composed by $\zz$ such that $B^{\downarrow}(\zz)\leq 0$ should have Gaussian fluctuations. This is related to the scaling relation $2 \chi=\xi$, where $\chi$ and $\xi$ are the longitudinal and transversal fluctuation exponents, respectively, which is conjectured to be universal (See, for instance, pg. 2089 in \cite{KZ} and pg. 589 in \cite{HoNe}.). For our exponential LPP model it is known that $\chi=1/3$ and $\xi=2/3$. Next we will see how use \eqref{eqBu} to prove the Gaussian fluctuations of the level sets of Busemann functions.    

Define 
$$\HH(t)=\{\zz\in\ZZ^2\,:\,B^{\downarrow}(\zz)\leq t\,\mbox{ and }\,B^{\downarrow}(\zz+\dd)>  t\}\,,$$ 
where $\dd=(1,1)$. We represent $\HH(t)$ as a down-right oriented bi-infinite path $\sigma_t=(\sigma_t(x))_{x\in\ZZ}$ in $\ZZ^2$, where we set $\sigma_t(0)$ to be the last point we have before passing through the diagonal. For instance, since $B(\0)=0$, we have that $\sigma_0(0)=-\ee_1$, $\sigma_0(1)=\0$ and $\sigma_0(2)=-\ee_2$. More generally we have the recursive relation: 
\begin{equation}\label{Bhoro}
\sigma_t(x+1)=\left\{\begin{array}{ll} \sigma_t(x)+\ee_1 & \mbox{ if }\, B^{\downarrow}( \sigma_t(x)+\ee_1)\leq t\,,\\
 \sigma_t(x)-\ee_2 & \mbox{ if }\, B^{\downarrow}( \sigma_t(x)+\ee_1)> t\,.\end{array}\right.
\end{equation}
The process $\sigma_t$ can be encoded as a particle process  $\zeta_t=(\zeta_t(x))_{x\in\ZZ}$ as follows:
\begin{equation}\label{Rost}
\zeta_t(x)=\left\{\begin{array}{ll} 0 & \mbox{ if }\, \sigma_t(x+1)-\sigma_t(x)=\ee_1\,,\\
1 & \mbox{ if }\, \sigma_t(x+1)-\sigma_t(x)=-\ee_2\,.\end{array}\right.
\end{equation}
We can think of $\zeta_t$ as a configuration of particles, where $\zeta_t(x) = 1$ means that there is a particle at site $x$ at time $t$, whereas $\zeta_t(x) = 0$ means that there is no particle at site $x$ at time $t$. This map between level sets and particle configuration is the so called Rost's correspondence \cite{Ro}. As a consequence of \eqref{eqBu} we get that: 
\begin{equation}\label{thm:Level} 
(\zeta_t)_{t\in\R}\stackrel{dist.}{=}(\eta^0_t)_{t\in\R}\,,
\end{equation}
where we recall that $(\eta^0_t)_{t\in\R}$ denotes the Palm version of the stationary totally asymmetric exclusion process with density $p=1/2$, conditioned to have a particle jumping from $0$ to $1$ at time zero. A straightforward consequence of \eqref{thm:Level} is that, after centering and rescaling in the standard way, the level set of the Busemann function will converge to a two-sided Brownian motion.

\begin{proof}[Proof of \eqref{thm:Level}] 
 Denote $\calR(\underline B^\downarrow)$ the deterministic map which associates the collection 
 $$\underline{B}^{\downarrow}:=\{ B^{\downarrow}\left(\zz\right)\,:\,\zz\in \ZZ^2\}\,$$
 to the particle process $(\zeta_t)_{t\in\R}$. It is straight forward to check that the map $\calR$ applied to $\underline{G}$ yields $(\eta^0_t)_{t\in\R}$:  
 $$\calR(\underline G)=(\eta^0_t)_{t\in\R}\,.$$
Thus, by \eqref{eqBu}, we have that $(\zeta_t)_{t\in\R}\stackrel{dist.}{=}(\eta^0_t)_{t\in\R}$.

\end{proof}

\subsection{Self-duality of the geodesic  tree}\label{SecDual}
By the coalescence property (iii), the collection of paths 
$$\cal L:=\{\gamma(\xx)\,:\,\xx\in\ZZ^2\}\,$$
is a.s. an up-right oriented tree, called the \emph{directional geodesic tree}. We also consider the collection of down-left oriented semi-infinite geodesics defined as  
$$\cal L^{\downarrow}:=\{\gamma^{\downarrow}(\xx)\,:\,\xx\in\ZZ^2\}\,.$$
It is clear that $\cal L$ and $\cal L^{\downarrow}$ have the same law, up to a rotation of $180$ degrees.

Let $\ZZ^{2*}$ denote the dual of $\ZZ^2$. We take as vertices the set $\{\zz^*=\zz+\frac{1}{2}\dd\,:\,\zz\in\ZZ^2\}$ (recall that $\dd=(1,1)$), and we join two such neighboring (distance $1$) vertices by a dual edge. Thus each edge of $\ZZ^2$ is bisected by a dual edge of $\ZZ^{2*}$, and vice-versa, which establishes a bijection (isomorphism) between edges and dual edges. Consider the last-passage percolation tree $\cal L$. The dual system $\cal L^*$ is defined as follows: in the case that an edge is in $\cal L$ then its dual is not in $\cal L^*$; in the case that an edge is not in $\cal L$ then its dual is in $\cal L^*$. Self-duality states that $\cal L$ and $\cal L^*$ have the same law, up to a rotation of $180$ degrees. The proof parallels the ideas in section 4.2 of \cite{FeMaPi}, where duality between geodesics and equilibrium competition interfaces was established, and relies on Lemma \ref{thm:BuBu}. 
\begin{lem}\label{thm:selfdual}
For the dual system, we have that 
$$\cal L^*\stackrel{dist.}{=}\cal L^{\downarrow}\,.$$
In particular, the dual system $\cal L^*$ is a.s. a tree and there is no bi-infinite maximizing path in $\cal L$.
\end{lem}

\begin{proof} For notational convenience, we will prove the equivalent statement that 
$$\cal L^{\downarrow*}\stackrel{dist.}{=}\cal L\,.$$
To prove that we first notice that the tree $\cal L^\downarrow$ is a deterministic function of the Busemann function $B^\downarrow$. In order to see this we use that
$$B^{\downarrow}(\xx)=\max\left\{ B^{\downarrow}(\xx-\ee_1), B^{\downarrow}(\xx-\ee_2)\right\}+ W_\xx\,.$$
Hence, for any down-left semi-infinite geodesic $\gamma^{\downarrow}(\xx)=(\xx_n)_{n\geq 0}$,  
\begin{equation}\label{backinfinite}
\xx_{n+1}=\argmax\left\{ B^{\downarrow}(\xx_n-\ee_1)\,,\, B^{\downarrow}(\xx_n-\ee_2)\right\}\,.
\end{equation}
(Notice that a similar property holds for finite geodesics.) By \eqref{backinfinite}, the tree $\cal L^{\downarrow}$ can be seen as the set composed of down-left oriented edges $(\xx,\ee_\xx)$ such that $\xx\in\ZZ^2$ and 
\begin{equation}\label{Bback}
\ee_{\xx}=\left\{\begin{array}{ll}\xx -\ee_1 & \mbox{ if }\, B^{\downarrow}(\xx-\ee_1)> B^{\downarrow}(\xx-\ee_2)\,,\\
\xx-\ee_2 & \mbox{ if }\,  B^{\downarrow}(\xx-\ee_2)> B^{\downarrow}(\xx-\ee_1)\,.\end{array}\right.
\end{equation}
Therefore, 
$$\cal L^{\downarrow}=\Upsilon( B^{\downarrow})\,$$
is a deterministic function $\Upsilon$ of $B^{\downarrow}=\{B^{\downarrow}(\xx)\,:\,\xx\in\ZZ^2\}$.

On the other hand, the dual system $\cal L^{\downarrow*}$ can be seen as the set composed of up-right oriented edges $(\xx^*,\ee_{\xx^*})$ such that
\begin{equation}\label{dualg}
\ee_{\xx^*}=\left\{\begin{array}{ll} \xx^*+\ee_1 & \mbox{ if }\, \ee_{\xx+\dd}=(\xx+\dd)-\ee_1\,,\\
\xx^*+\ee_2 & \mbox{ if }\, \ee_{\xx+\dd}=(\xx+\dd)-\ee_2\,.\end{array}\right.
\end{equation}
In other words, the edge in $\cal L^{\downarrow*}$ starting at vertex $\xx^*=\xx+\frac{1}{2}\dd$ will point up or right if the edge in $\cal L^{\downarrow}$ starting at $\xx+\dd$ points down or left, respectively. Now, by \eqref{Bback} and \eqref{dualg},
$$\ee_{\xx^*}=\left\{\begin{array}{ll} \xx^*+\ee_1 & \mbox{ if }\,  B^{\downarrow*}(\xx^*+\ee_1)< B^{\downarrow*}(\xx^*+\ee_2)\,,\\
\xx^*+\ee_2 & \mbox{ if }\, B^{\downarrow*}(\xx^*+\ee_2)< B^{\downarrow*}(\xx^*+\ee_1)\,,\end{array}\right.$$
where $B^{\downarrow*}(\xx^*):= B^{\downarrow}(\xx)$. Let $\phi:\xx\in\ZZ^2\mapsto\phi(\xx):=(-\xx)^*\in\ZZ^{2*}$ and set 
$$\tilde B(\xx):=- B^{\downarrow*}(\phi(\xx))\,.$$
Then we have that $\phi^{-1}(\cal L^{\downarrow*})$ can be represented as the set composed of down-left oriented edges $(\xx,\ee_{\xx})$ such that
\begin{equation*}
\ee_{\xx}=\left\{\begin{array}{ll}\xx -\ee_1 & \mbox{ if }\, \tilde B(\xx-\ee_1)>\tilde B(\xx-\ee_2)\,,\\
\xx-\ee_2 & \mbox{ if }\,  \tilde B(\xx-\ee_2)>\tilde B(\xx-\ee_1)\,.\end{array}\right.
\end{equation*}
Or, equivalently, 
$$\phi^{-1}(\cal L^{\downarrow*})=\Upsilon(\tilde B)\,.$$
By Lemma \ref{thm:BuBu},
\begin{equation*}
\{\tilde B(\xx)\,:\,\xx\in\ZZ^2\}\stackrel{dist.}{=}\{B^{\downarrow}(\xx)\,:\,\xx\in\ZZ^2\}\,.
\end{equation*}
Hence,
$$\phi^{-1}(\cal L^{\downarrow*})=\Upsilon(\tilde B)\stackrel{dist.}{=}\Upsilon(B^{\downarrow})=\cal L^{\downarrow}\,,$$
and the proof of self-duality is completed.

By self-duality, all almost sure statements for $\cal L$ also hold for $\cal L^*$. Therefore, a.s. $\cal L^*$ is a tree. If, with positive probability, there were a bi-infinite path in $\cal L$, then the dual system $\cal L^*$ would be split into two disjoint parts, which can not happen, since $\cal L^*$ is a.s. a tree.    

\end{proof}

\subsection{Proof of the duality formula}
The equivalence between the stationary LPP model with boundary and Busemann functions allows us to interpret exit points as \emph{crossing points} of semi-infinite geodesics. For $x,n\geq 1$, let $Z^{\downarrow}(x,n)$ denote the first point in $\gamma^{\downarrow}((x,n))$ (following the down-left orientation) that intersects $[1,x]\times\{0\}\cup\{0\}\times[1,n]$. Notice that this intersection has to be transversal to the axis. Again, to distinguish between crossings via the horizontal or the vertical axis, we introduce a non-zero integer-valued random variable $Z^{\downarrow}$ such that if $Z^{\downarrow}>0$ then the crossing point is $(Z^{\downarrow},0)$, while if $Z^{\downarrow}<0$ then the crossing point is $(0,-Z^{\downarrow})$. Define the \emph{crossing-point process} as
\begin{equation*}%\label{defpoint}
\cal Z^{\downarrow}_n:=(\zeta^{\downarrow}_n(z)\,,z\in[-n,\infty))\in\{0,1\}^{[-n,\infty)}\,,
\end{equation*}  
where, for fixed $n\geq 1$,
$$\zeta^{\downarrow}_n(z)=\left\{\begin{array}{ll}1 & \mbox{ if $z=Z^{\downarrow}(x,n)$ for some $x\in [1,\infty)$}\,,\\ 
0 & \mbox{otherwise}\,.\end{array}\right.$$ 
In \cite{CaPi}, it was proved that
\begin{equation}\label{maxcross}
\cal Z_n^{\downarrow}\stackrel{dist.}{=}\cal Z_n\,.
\end{equation}

A key observation is that the coalescence time $T^{\downarrow*}_m:=T_m(\cal L^{\downarrow*})$ of the dual tree and the crossing point process $\cal Z_n^{\downarrow}$ are related by: 
\begin{equation}\label{noncross}
\{T^{\downarrow*}_m< n\}=\{\cal Z^{\downarrow}_n([-m,m])=0\}\,.
\end{equation}
This is a topological consequence of the fact that, by definition, $\cal L^{\downarrow}$ and $\cal L^{\downarrow*}$ do not cross each other. Hence, if $T^{\downarrow*}_m<n$, then the dual paths emanating from $(\mm^h)^*$ and $(\mm^v)^*$ prevent that $Z^{\downarrow}(x,n)\in[-m,m]$ for any $x \geq 1$, and vice-versa (recall that $Z^\downarrow$ is the transversal intersection point). Now we are able to prove the duality formula.

\begin{proof}[Proof of \eqref{dualfor}]
Recall that $T_m=T_m(\cal L)$ (it is a deterministic function of the tree $\cal L$). By Lemma \ref{thm:selfdual}, we have that 
\begin{equation}\label{dualcoal}
T_m(\cal L)\stackrel{dist.}{=}T_m(\cal L^{\downarrow*})\,.
\end{equation}
Therefore, by \eqref{dualcoal}, \eqref{noncross} and \eqref{maxcross} (in this order),
\begin{eqnarray*}
\P\left(T_m< n\right)&=&\P\left(T^{\downarrow*}_m< n\right)\\
&=&\P\left(\cal Z^{\downarrow}_n[-m,m]=0\right)\\
&=&\P\left(\cal Z_n[-m,m]=0\right)\,,
\end{eqnarray*}
and the proof of \eqref{dualfor} is finished.

\end{proof}

\subsection{Proof of the lower bounds}

\begin{proof}[Proof of \eqref{coaltail}] Denote $Z_n:=Z(n,n)$. By Theorem 2.2 in \cite{BaCaSe} \footnote{Notice that $|Z_n|=Z_{n+}+Z_{n-}$ and that $Z_{n+}\stackrel{dist.}{=}Z_{n-}$.} there exists a constant $c>0$ such that if 
$$\lim_{m\to\infty}n/m^{3/2}=r>0$$
(where $n=n(m)$), then
$$\limsup_{m\to\infty}\P\left(|Z_n|\geq m \right)\leq cr^{2}\,.$$
On the other hand, 
$$\P\left(\cal Z_{n+1}[-m,m]\geq 1\right)\geq \P\left(Z_{n+1}\in [-m,m]\right)\,.$$
Together with the duality formula, this yields  
\begin{equation*}
\P\left(T_m> n\right)\geq \P\left(Z_{n+1}\in [-m,m]\right)\,,
\end{equation*}
and hence
$$\liminf_{m\to\infty}\P\left(T_m> n\right)\geq 1-\limsup_{m\to\infty}\P\left(|Z_{n+1}| >m\right)\geq 1-cr^2$$
as soon as $n/m^{3/2}\to r$. 

\end{proof}

The last-passage time $\bar L$ has  a variational representation given by    
\begin{equation}\label{lastvar}
\bar L(x,n)=\max_{z\in[-n,x]}\left\{M(z)+L_z(x,n)\right\}\,,\mbox{ for $x,n\geq 1$}\,,
\end{equation}
where $M(z)$ is the sum of the (i.i.d. $\Exp(1/2)$) passage times along the boundary, 
$$M(z):=\left\{\begin{array}{ll} 0\,, & \mbox{ if }\, z =0\,;\\
\sum_{k=1}^z\bar W_{(k,0)}\,, & \mbox{ if }\, z > 0\,;\\
\sum_{k=1}^{-z}\bar W_{(0,k)}\,,& \mbox{ if }\,z<0\,,\end{array}\right.$$
and 
$$L_z(x,n):=\left\{\begin{array}{ll} L(\0,(x,n))\,, & \mbox{ if }\, z =0\,;\\
L((z,0),(x,n))\,, & \mbox{ if }\, z > 0\,;\\
L((0,-z),(x,n))\,,& \mbox{ if }\,z<0\,.\end{array}\right.$$
Therefore
$$\bar L(x,n)=M(Z)+L_{Z}(x,n)\,,$$
or, in other words, exit-points of geodesics are locations of maxima: 
\begin{equation}\label{argvar}
Z(x,n)=\argmax_{z\in [-n,x]}\left\{M(z)+L_z(x,n)\right\}\,.
\end{equation}
This variational representation for exit points, together with the scaling limit of last-passage times, implies a limit theorem for $Z_n$. 

\begin{lem}\label{unique}
Let $ \cal B$ be a two-sided standard Brownian motion and let $\cal A$ be an independent Airy$_2$ process. Then a.s. there is a unique location $U\in\R$ such that  
$$U:=\argmax_{u\in\R}\left\{\sqrt{2} \cal B(u) + \cal A(u)-u^2\right\}\,,$$
\end{lem}

\begin{proof}
We apply the method of proof developed in \cite{Pi} to show uniqueness of the location of maxima for a continuous process. There, it was proven that uniqueness of the location of the maxima of a continuous  process $X$ is equivalent to the existence of the derivative of the function 
$$m(a):= \E \left(\max_{u\in[0,t]}\left\{X(u)+au\right\}\right)\,,\,\mbox{ for }a\in\R\,,$$   
at $a=0$. This result was also generalize for $X(u)=\cal B(u)-u^2$ and $X(u)=\cal A(u)-u^2$, where the maximisation was taken over $u\in\R$. We use the same idea of proof for $X(u)=\sqrt{2}\cal B(u)+\cal A(u)-u^2$. Indeed, by completing the square, we get that
\begin{equation}\label{square}
\sqrt{2}\cal B(u)+\cal A(u)-u^2+au=\sqrt{2} \cal B(u)+\cal A(u)-\left(u-\frac{a}{2}\right)^{2}+\frac{a^2}{4}\,.
\end{equation}
We note that 
$$\cal B(s+a/2)-\cal B(a/2)\stackrel{dist.}{=}\cal B(s)$$ 
(by shift invariance), and that 
$$\cal A(s+a/2)\stackrel{dist.}{=}\cal A(s)\,$$
(by stationarity). 
Set $s=u-a/2$ and add and subtract $\sqrt{2}\cal B(a/2)$ to \eqref{square}. Then these distributional invariances imply that 
$$m(a)=m(0)+\frac{a^2}{4}\,,$$ 
which shows differentiability at $a=0$ and, as a consequence, a.s. uniqueness of the location of the maxima.

\end{proof}

\begin{lem}\label{thm:scalargmax}
Define
$$U_n:=\frac{Z_n}{2^{5/3}n^{2/3}}\,,$$
Then
$$\lim_{n\to\infty}U_n\stackrel{dist.}{=}U\,.$$  
\end{lem}

\begin{proof} 
We present a sketch of the proof and leave further details to the reader. It follows a similar structure as in the proof of convergence of the location of maxima in the point to line LPP model, developed in \cite{Jo}. The first ingredient is the following functional limit result   
\begin{equation}\label{lppairy}
\lim_{n\to\infty}\cal A_n(u)\stackrel{dist.}{=}\cal A(u)\,,
\end{equation}
where
$$\cal A_n(u):=\frac{L_{2^{5/3}un^{2/3}}(n,n)-\left(4n-2^{8/3}un^{2/3}\right)+2^{4/3}u^2 n^{1/3}}{2^{4/3} n^{1/3}}\,.$$
For finite dimensional convergence see \cite{Jo}, and for tightness see \cite{CaPi1}. By the functional central limit theorem, we have that
\begin{equation}\label{lppbrow}
\lim_{n\to\infty}\cal B_n(u)\stackrel{dist.}{=}\sqrt{2}\cal B(u)\,,
\end{equation}
where
$$\cal B_n(u):=\frac{M(2^{5/3}un^{2/3})-2^{8/3}un^{2/3}}{2^{4/3}n^{1/3}}\,.$$
Let
$$\cal C_n:=\frac{\bar L(n,n)-4n}{2^{4/3}n^{1/3}}\,.$$
By \eqref{lastvar}, we have that (for $c=2^{-5/3}$)
$$\cal C_n=\max_{u\leq cn^{1/3}}\left\{\cal B_n(u)+\cal A_n(u)-u^2\right\}\,,$$
and hence (notice that $\cal A_n$ and $\cal B_n$ are independent), 
\begin{equation}\label{lppa+b}
\lim_{n\to\infty}\cal C_n\stackrel{dist.}{=}\max_{u\in\R}\left\{\sqrt{2}\cal B(u)+\cal A(u)-u^2\right\}\,.
\end{equation}
See \cite{BaFePe} for a description of the limit law of $\cal C_n$, and \cite{QuaRe} for more details on variational problems involving the Airy$_2$ process and the Brownian motion.

By Theorem 2.2 in \cite{BaCaSe}, $(U_n)_{n\geq 1}$ is tight and, by \eqref{argvar}, 
\begin{equation}\label{scalargvar}
U_n=\arg\max_{u\leq cn^{1/3}}\left\{\cal B_n(u)+\cal A_n(u)-u^2\right\}\,.
\end{equation}
Therefore, Theorem \ref{thm:scalargmax} will follow as soon as the location of maxima of the limit process is a.s. unique (to have continuity of the $\arg\max$ functional), which is given by Lemma \ref{unique}. 

\end{proof}

Now, we apply Lemma \ref{thm:scalargmax} to lower bound $\G$.

\begin{proof}[Proof of \eqref{lowtail}] 
As we saw in the proof of the previous corollary,
$$\P\left(T_m> n\right)\geq \P\left(Z_{n+1}\in (-m,m]\right)\,,$$
and hence
$$\P\left(\frac{T_m}{2^{-5/2}m^{3/2}}> \frac{n}{2^{-5/2}m^{3/2}}\right)\geq\P\left(U_{n+1}\in \left(-\frac{m}{2^{5/3}n^{2/3}},\frac{m}{2^{5/3}n^{2/3}}\right]\right)\,.$$
If we take $m,n$ such that $n/2^{-5/2}m^{3/2}\to r$, then $m/2^{5/3}n^{2/3}\to r^{-2/3}$. Thus, by Lemma \ref{thm:scalargmax},     
$$\G(r)\geq\P\left(U\in(-r^{-2/3},r^{-2/3}] \right)=\F(r^{-2/3})-\F(-r^{-2/3})\,.$$

\end{proof} 

\section{Final Comments}\label{sec:final}
\subsection{Upper bounds}
To get sharp upper bounds for coalescence times one needs to show that 
\begin{equation}\label{upper}
\limsup_{n\to\infty}\P\left(\calZ_n([-\delta n^{2/3},\delta n^{2/3}])\geq 1\right)\leq c\delta\,,
\end{equation}
for some fixed constant $c>0$ and small enough $\delta>0$. One possible approach is to parallel the arguments developed in \cite{CaGro} to bound the probability that $Z_n\in [0,\delta n^{2/3}]$. However, an extra (and non trivial) effort will be necessary since one will need uniform control over the whole exit point process (not only at single location).    

\subsection{Duality in the scaling limit}
Define the rescaled processes   
$$\cal A_n(u,v):=\frac{L_{2^{5/3}un^{2/3}}(n+2^{5/3}vn^{2/3},n)-\left(4n+2^{8/3}(v-u)n^{2/3}\right)+2^{4/3}(u-v)^2 n^{1/3}}{2^{4/3} n^{1/3}}\,,$$
and
$$\cal C_n(v):=\frac{\bar L(n,n+2^{5/3}vn^{2/3})-(4n+2^{8/3}vn^{2/3})}{2^{4/3}n^{1/3}}\,.$$
By \eqref{lastvar}, we have that 
$$\cal C_n(v)=\max_{u\leq n^{1/3}}\left\{\cal B_n(u)+\cal A_n(u,v)-(u-v)^2\right\}\,.$$
The process $\cal C_n(v)$ has a limit \cite{BaFePe}
$$\lim_{n\to\infty}\cal C_n(v)\stackrel{dist.}{=}\cal C(v)\,,$$
whose finite dimensional distributions are also expressed in terms of Fredholm determinants. It is known that the sequence $(\cal A_n)_{n\geq 1}$ is tight (in the space of two parameter continuous processes), although no rigorous result on the convergence of finite dimensional distributions is available \cite{CaPi1, CoQua}. For fixed $u$, it is not hard to see that, for fixed $v\in\R$,  
$$\lim_{n\to\infty}\cal A_n(u,v)\stackrel{dist.}{=}\cal A(u-v)\,,$$
as a process in $u\in\R$. It is conjectured that $\cal A_n(u,v)$ indeed converges to a two parameter process $(\cal A(u,v)\,, (u,v)\in\R^2)$, called the \emph{Airy$_2$ sheet} \cite{CoQua}. The Airy$_2$ sheet is symmetric and  stationary process with continuous paths. These limit processes are related to each other by the variational relation   
$$\cal C(v)-\cal C(0)\stackrel{dist.}{=}\sqrt{2}\cal B(v)\,,\mbox{ for }v\in\R\,\mbox{ (as process) }$$
where 
$$\cal C(v):=\max_{u\in\R}\left\{\sqrt{2}\cal B(u)+\cal A(u,v)-(u-v)^2\right\}\,.$$

Consider the jump process $(U(v)\,,\,v\in\R)$ which runs through the (right-most) locations of maxima: 
$$U(v):=\sup\arg\max_{u\in\R}\left\{\sqrt{2}\cal B(u)+\cal A(u,v)-(u-v)^2\right\}\,.$$
By stationarity of the Airy$_2$ sheet and shift invariance of the two-sided Brownian motion, the process $(U(v)-v\,,\,v\in\R)$ will be stationary. (It is also known that $\E U(0)=0$, and hence $\E U(v)=v$.) Define the counting process $\cal U:=(\zeta(u)\,,\,u\in \R)$ induced by the locations of maxima: 
$$\zeta(u)=\left\{\begin{array}{ll}1 & \mbox{ if $u=U(v)$ for some $v\in\R$}\,,\\ 0 & \mbox{otherwise}\,.\end{array}\right.
$$
and 
$$\cal U(A):=\sum_{u\in A}\zeta(u)\,.$$ 

Based on the variational representation \eqref{argvar} of exit points, we conjecture that the exit-point counting process $\cal Z_n$, rescaled by $2^{5/3}n^{2/3}$, converges to $\cal U$. By duality \eqref{dualfor}, if this conjecture is true, one gets the existence of the limiting distributuion, 
$$T\stackrel{dist.}{=}\lim_{m\to\infty}\frac{T_m}{2^{-5/2}m^{3/2}}\,,$$
and that
\begin{equation}\label{eq:dualimit}
\P\left(T\leq r\right)=\P\left(\cal U\left((-r^{-2/3},r^{-2/3}]\right)=0\right)\,.
\end{equation}
We also expect that 
\begin{equation*}
\lim_{r\to\infty}\frac{\G( r)}{r^{2/3}}= 2\lambda\,. 
\end{equation*}
where
$$\lambda:=\lim_{\delta\to 0^+}\frac{\P\left(\cal U\left((0,\delta]\right)\geq 1\right)}{\delta}\,.$$

%A natural guess for $\lambda$ is $\lambda=f(0)$, where $f$ is the density function of $U$. We note that $U(0)=\delta V_\delta$ where
%\begin{eqnarray*}
%V_\delta&:=&\arg\max_{v\in\R}\left\{\sqrt{2}\cal B(\delta v)+\cal A(\delta v)-(\delta v)^2\right\}\\
%&=&\arg\max_{v\in\R}\left\{\sqrt{2}\delta^{-1/2}\cal B(\delta v)+\delta^{-1/2}(\cal A(\delta v)-A(0))-\delta^{3/2} u^2\right\}\\
%&\stackrel{dist.}{\to}&\arg\max_{v\in\R}\left\{2\cal B(v)-\delta^{3/2} u^2\right\}\,\,\,\mbox{(as $\delta\to 0$)}\\
%&\stackrel{dist.}{=}&\delta^{-1}\arg\max_{v\in\R}\left\{2\cal B(v)-v^2\right\}=:\delta^{-1}\bar V\,.
%\end{eqnarray*}
%The convergence in the third equation needs to be formally justified, but since the Airy process behaves locally as Brownian motion \cite{Ha,CaPi1}, it seems to be a reasonable guess at a heuristic level. If this is true, we have that $f(0)$ equals $\bar f(0)$, where $\bar f$ is the density function of $\bar V$. The random variable $2^{-2/3}\bar{V}$ has the so-called Chernoff distribution \cite{GroWe}. We also wonder if the same argument (local convergence) indicates that $\lambda=2^{2/3}\bar \lambda$, where $\bar\lambda$ is the intensity of the point process composed by the locations of a Brownian motion minus a parabola \cite{Gro}.       

\subsection{The polymer web} The Airy$_2$ sheet can also be seen as a space-time parameter process $\cal A(s,u;t,v)$, where $\cal A(u,v)=\cal A(0,u;1,v)$. This space-time process is conjectured to be the space-time scaling limit of last-passage percolation models, and also of solutions to the Kadar-Parisi-Zhang equation \cite{CoQua}. It induces a random semi-group $T_{s,t}$, acting on functions $f$ by the variational formula
$$\cal C_{s,t}(f)(v):=\max_{u\in\R}\left\{f(u)+\cal A(s,u;t,v) -\frac{(u-v)^2}{(t-s)}\right\}\,.$$
The two-sided Brownian motion is a fixed point in the sense that 
$$\cal C_{0,t}(\cal B)(v)-\cal C_{0,t}(\cal B)(0)\stackrel{dist.}{=}\cal B(v)\mbox{ for all $t\geq 0$}\,.$$  
In this context, one could consider the time process composed by counting measures $(\cal U_t\,,t\geq 0)$ induced by Dirac deltas located at maxima of $\cal C_{0,t}(\cal B)$:
\begin{equation}\label{polyweb}
U(v,t):=\sup\arg\max_{u\in\R}\left\{\sqrt{2}\cal B(u)+\cal A(0,u;t,v)-\frac{(u-v)^2}{t}\right\}\,
\end{equation}
The parabolic term forces $U(v,t)$ to be close to $v$, and its effect has a decreasing influence as $t\to\infty$, which implies that the locations will became more and more sparse as $t\to\infty$. Therefore, it is natural to think in terms of the trajectories of the locations, and that these locations will coalesce as time passes. We conjectured that this collection of coalescing trajectories, which we call the polymer web, is the scaling limit of the directional geodesic tree.

\subsection{Coalescence times and local equilibrium}
The LPP model can be represented as a discrete time Markov interacting system $(M_n\,:\,n\geq 0)$ on $[0,\infty)^\ZZ$ \cite{CaPi}. (See also \cite{AlDi} for the Hammersley LPP model and its particle system interpretation.) At time zero we start with a collection of non-negative weights $\{W_i\,:\,i\in \ZZ\}$. We define the weight (or mass) of the interval $(a,b]$ at time zero as 
$$M(a,b]=M_0(a,b]:=\sum_{i=a+1}^bW_i\,.$$
At time $n\geq 1$, we define the weight of the interval $(a,b]$ as 
$$M_n(a,b]:=\bar L_{M}(b,n)-\bar L_M(a,n)\,$$
where \footnote{We assume that $\liminf_{k\to\infty}\frac{\sum_{i=-1}^{-k}W_i}{k}>1$ so that the maxima is indeed attained on a compact set. }
$$\bar L_M(x,n):=\max_{z\leq x}\left\{M(z)+L_z(x,n)\right\}\,,$$
and
$$M(z):=\left\{\begin{array}{ll} 0\,, & \mbox{ if }\, z =0\,;\\
\sum_{k=1}^z W_k\,, & \mbox{ if }\, z > 0\,;\\
-\sum_{k=1}^{-z} W_{k}\,,& \mbox{ if }\,z<0\,.\end{array}\right.$$
The last-passage time $L$ can be recovered by choosing an initial weight configuration with infinite mass at negative  sites, and with i.i.d exponential weights of parameter one at non-negative sites.  

Time stationary ergodic measures on $[0,\infty)^\ZZ$ for this Markov system are represented by i.i.d. collections of exponentials random weights of parameter $\rho\in(0,1)$: if $\{W_i\,:\,i\in \ZZ\}$ is distributed according to an i.i.d. collection of $\Exp(\rho)$ random variables, then 
$$M_n  \stackrel{dist.}{=} M_0\,,\mbox{ for all }n\geq 0\,.$$  
Local equilibrium of the LPP interacting system is described by Busemann functions \cite{CaPi}: 
$$\lim_{n\to\infty}L(\0,(n+h,n))-L(\0,(n,n))\stackrel{dist.}{=}B^{\downarrow}(\0,(0,h))\stackrel{dist.}{=}M_0(0,h]\,,$$ 
with $\rho=1/2$. If one moves the origin to $-\nn:=-(n,n)$ (and starts the system at time $-n$) then the convergence becomes a.s.: 
$$\lim_{n\to\infty}L(-\nn,(h,0))-L(-\nn,\0)\stackrel{a.s.}{=}L(\cc_h,(0,h))-L(\cc_h,\0)=B^{\downarrow}(\0,(0,h))\,,$$ 
where $\cc^{\downarrow}_h:=\cc^{\downarrow}(\0,(0,h))$. Thus, the coalescence time also describes how far in the past one needs to start the process to see local equilibrium in the present. In this sense, it would be interesting to analyze coalescence times and duality \eqref{dualfor} in the framework of  relaxation and mixing times for Markov processes.

\end{document}